\documentclass[final]{siamltex}

\usepackage{times,amscd,amsmath,amssymb,amsbsy,epsfig}

\DeclareFontFamily{OML}{script}{}
\DeclareFontShape{OML}{script}{m}{it}
{ <5-20> rsfs10 }{}
\DeclareMathAlphabet{\mathscript}{OML}{script}{m}{it}

\newcommand{\re}[1]{\mbox{$($\ref{#1}$)$}}

\newtheorem{thm}{Theorem}[section]
\newtheorem{cor}[thm]{Corollary}
\newtheorem{lma}[thm]{Lemma}
\newcommand{\pf}{\par\noindent {\textit {Proof. }}}
\newcommand{\q}{~\mbox{\hspace{2em}\frame{\rule{0ex}{1.5ex}\mbox{\hspace{1ex}}}}\xx}
\numberwithin{equation}{section}
\DeclareMathOperator{\esssup}{ess\,sup}

\newcommand{\xx}{\vspace*{2ex}}
\newcommand{\m}{\hspace{1em}}
\newcommand{\mm}{\hspace{2em}}
\newcommand{\on}{\quad\text{on}\;\;}
\newcommand{\iin}{\quad\text{in}\;\;}

\newcommand{\at}{\quad\mbox{at}\;\;}
\newcommand{\as}{\quad\text{as}\;\;}
\newcommand{\iif}{\quad\text{if}\;\;}
\renewcommand{\le}{\leqslant}
\newcommand\beq{\begin{equation}}
\newcommand\eeq{\end{equation}}
\newcommand\be{\begin{equation}}
\newcommand\ee{\end{equation}}
\newcommand\beqn{\begin{equation*}}
\newcommand\eeqn{\end{equation*}}
\newcommand\bea{\begin{eqnarray}}
\newcommand\eea{\end{eqnarray}}
\newcommand\beaa{\begin{eqnarray*}}
\newcommand\eeaa{\end{eqnarray*}}

\newcommand{\cF}{\mathcal{F}}\newcommand{\sF}{\mathscript{F}}

\newcommand{\sL}{\mathscript{L}}

\newcommand{\bv}{\mathbf{v}}

\newcommand{\gm}{\gamma}

\newcommand{\tdR}{\tilde{R}}\newcommand{\tdrho}{\tilde{\rho}}
\newcommand{\tdM}{\tilde{M}}\newcommand{\tdK}{\tilde{K}}

\newcommand{\ld}{\lambda_}
\newcommand{\pd}[2]{\dfrac{\partial#1}{\partial#2}}
\newcommand{\dd}[2]{\dfrac{d#1}{d#2}}

\newcommand{\divp}[1]{\frac{1}{r}\frac{\partial}{\partial r}\left(r #1 \right)}

\newcommand{\eps}{\varepsilon}

\newcommand{\picturesAB}[3]{
\centerline{\raise #3 \hbox{(a)}
\psfig{file=#1,height=#3}
\hspace*{.2in}
\raise #3 \hbox{(b)}
\psfig{file=#2,height=#3}
}}
\newcommand{\picturesABC}[4]{
\centerline{\raise #4 \hbox{(a)}
\psfig{file=#1,height=#4}
\hspace*{.2in}
\raise #4 \hbox{(b)}
\psfig{file=#2,height=#4}
\hspace*{.2in}
\raise #4\hbox{(c)}
\psfig{file=#3,height=#4}
}}

\newcommand{\kineticrhon}{\frac{k_{\rho} w }{w+K_{w\rho}} f (1-\frac{\rho}{\rho_m}) - \lambda_{\rho} \rho }
\newcommand{\kineticwn}{k_{w} b \big( (1-\gamma)w_b - w \big) - \left[ \big(\ld{wf} f+ \ld{wm} m\big)  \left(1 + \dfrac{\ld{ww} p}{1+ p}\right) + \ld{wm} \right] w   }
\newcommand{\kineticpn}{k_p m  G_p(w) - \dfrac{ \ld{pf}f p}{1 + p} -\ld{p} p}
\newcommand{\kineticen}{k_e m G_e(w) - (\ld{en}n + \ld{eb}b +\ld{e}) e}
\newcommand{\kineticmn}{\frac{ k_{m}   b p}{1+p} - \ld m m \left(1 + \lambda_d D(w) \right)}
\newcommand{\kineticfn}{k_f G_f(w)   f \left(1-\frac{f}{f_m}\right) - \ld{f} f (1 + \lambda_d D(w) )}
\newcommand{\kineticnn}{(k_{nb}  b + k_{n}  n) \frac{e}{1 + e} - (\ld{nb} b + \ld{nn} n)n}
\newcommand{\kineticbn}{k_{b} G_b(w ) b (1- b ) +  G_b(w) (\ld{nb} b + \ld{nn} n)n}
\newcommand{\conv}[1]{\frac{1}{r}\pd{}{r} \big(r #1 v\big)}
\newcommand{\bddtaxis}[2]{\frac{#1 {\partial #2}/{\partial r}}{\sqrt{1+ k_{sg}\left|{\partial #2}/{\partial r}\right|^2} } }
\newcommand{\bddtaxist}[2]{\frac{#1 {\partial #2}/{\partial \xi}}{\sqrt{1+ k_{sg}(t)  \left|{\partial#2}/{\partial \xi}\right|^2 } } }
\newcommand{\diffct}[1]{D_{#1}(t)}
\newcommand{\chit}[1]{\chi_{#1}(t)}

\title{Analysis of a mathematical model of \\ischemic cutaneous wounds}

\author{Avner Friedman\thanks{Mathematical Biosciences Institute,
and Department of Mathematics, Ohio State University, Columbus, Ohio 43210
(afriedman@math.ohio-state.edu)}\and Bei Hu\thanks{Department of Mathematics,
University of Notre Dame, Notre Dame, Indiana 46556 (b1hu@nd.edu) }\and
Chuan Xue\thanks{Mathematical Biosciences Institute, Ohio State University,
Columbus, Ohio 43210 (cxue@mbi.osu.edu)} }

\begin{document}

\maketitle

\begin{abstract}
Chronic wounds represent a major public health problem affecting 6.5 million
people in the United States. Ischemia represents a serious complicating factor
in wound healing. In this paper we analyze a recently developed mathematical
model of ischemic dermal wounds. The model consists of a coupled system of
partial differential equations in the
partially healed region, with the wound boundary as a free boundary.
The extracellular matrix (ECM) is assumed to be viscoelastic, and the
free boundary moves with the velocity of the ECM at the boundary of the open
wound. The model equations involve the concentrations of oxygen, cytokines, and the densities
of several types of cells. The ischemic level is represented by a parameter which appears
in the boundary conditions, $0\leq \gamma < 1$; $\gamma$ near 1 corresponds to extreme ischemia and
$\gamma=0$ corresponds to normal non-ischemic conditions. We establish global existence and uniqueness of
the free boundary problem and study the dependence of the free boundary on $\gamma$.
\end{abstract}



\section{Introduction}

Wound healing represents the outcome of
a large number of interrelated biological events that are orchestrated
over a temporal sequence in response to injury and its microenvironment.
The process involves interactions among different soluble chemical
mediators, different types of cells, and the extracellular matrix (ECM).
Among the various factors that affect the healing of a wound, the tissue oxygen
level is a key determinant \cite{Gordillo:2003:RER,Sen:2009:WHE}. Although hypoxia
is generally recognized as a physiological cue to induce angiogenesis
\cite{Colgan:2007:HLD,Semenza:2000:HUT,Pugh:2003:RAH,Liao:2007:HKR}, severe hypoxia
cannot sustain the growth of functional blood vessels
\cite{Hopf:2005:HA, Allen:1997:WHN, Gibson:1997:IOT, Oberringer:2005:DCC, Safran:2003:HHM}.

There have been several mathematical models of wound healing which incorporated the effect
of angiogenesis \cite{Pettet:1996:MWA,Pettet:1996:RAW,Byrne:2000:MMA,Schugart:2008:WAF}.
Mathematical models of angiogenic networks, such as through the induction of vascular
networks by vascular endothelial growth factors (VEGFs) \cite{Dor:2003:IVN,Dor:2003:MVN},
were developed by McDougall and coworkers \cite{McDougall:2002:MMF,Stephanou:2005:MMF},
based in part on the work of Anderson and Chaplain \cite{Anderson:1998:CDM}, in connection with chemotherapeutic
strategies. The role of oxygen in wound healing was explicitly incorporated in the works of Byrne et al. \cite{Byrne:2000:MMA}
and Schugart et al. \cite{Schugart:2008:WAF}. In particular, it was demonstrated in \cite{Schugart:2008:WAF}
that enhanced healing can be achieved by moderate hyperoxic treatments.
In \cite{Roy:2009:CPM}, the impairment of dermal wound healing due to ischemic conditions was
addressed in a pre-clinical experimental model. In a more recent work \cite{Xue:2009:MMI},
Xue, Friedman and Sen developed a mathematical model of ischemic dermal wound-healing. The model consists of a
system of PDEs in the partially healed region which is modeled as a viscoelastic medium with a free boundary
surrounding the open wound. Simulations of the model were shown to be in  agreement with the experimental results in \cite{Roy:2009:CPM}.

In this paper we study the model in \cite{Xue:2009:MMI} by mathematical analysis. In particular we prove that the free
boundary problem developed in that model has a unique global solution, and that the open wound does not close under extreme
ischemic conditions. We also show, by simulations, that non-ischemic wounds do heal. In Section \ref{2} we formulate the mathematical model
for a radially symmetric geometry as in \cite{Xue:2009:MMI}.
The ischemic level is determined by a parameter $\gamma$, $0\leq \gm \leq 1$;
$\gamma$ near 1 corresponds to extreme ischemia and $\gamma=0$ corresponds to normal non-ischemic conditions.
In Section \ref{3} we show that the free boundary is monotone decreasing, and in Section \ref{4} we derive \textit{a priori} estimates. In
Section \ref{5} we transform the free boundary problem into a problem in a fixed domain; this is a convenient form for proving, in Section
\ref{6}, local existence and uniqueness of a solution. The extension of the solution to all $t>0$ is
also established in Section \ref{6}, by using the \textit{a priori}
estimates derived in Section \ref{4}. In Section \ref{7} we consider
the case of extreme ischemia (namely, $\gamma$ near 1) and prove
that the wound's boundary stops decreasing after some finite time.
In Section \ref{8} we establish some properties of the solution for
wounds that do not heal. Section \ref{9} simulates the radius of the
wound when the parameter of the system are chosen, as in
\cite{Xue:2009:MMI}, based on biological literature. The simulations
suggest the following conjecture: there exists a parameter
$\gm^{\ast}$ such that wounds heal if $0\leq \gm < \gm^{\ast}$ and do
not heal if $\gm^{\ast}< \gm \leq 1$.

\section{The mathematical model}\label{2}

It is assumed that the dermal tissue is in a circular domain $\{(r,\theta; \; r\leq L)\}$ and
the open wound at time $t$ is a disc $\{(r,\theta; \; r<R(t))\}$ with initial
radius $R(0)<L$. The partially healed tissue is the annulus $\Omega(t) = \{(r,\theta; \; R(t)\leq r\leq L)\}$.
We introduce the following variables:
\begin{itemize}
\item
Chemicals:
\begin{itemize}
\item[] $w(r,t)$: concentration of tissue oxygen
\item[] $e(r,t)$: concentration of Vascular Endothelial Growth Factor (VEGF)
\item[] $p(r,t)$: concentration of Platelet Derived Growth Factor (PDGF)
\end{itemize}
\item Cells, blood vessels and matrix
\begin{itemize}
\item[] $m(r,t)$: density of macrophages
\item[] $f(r,t)$: density of fibroblasts
\item[] $n(r,t)$: density of capillary tips
\item[] $b(r,t)$: density of capillary sprouts
\item[] $\rho(r,t)$: density of the ECM 
\item[] $v(r,t)$: velocity of the ECM
\end{itemize}
\end{itemize}
In homeostasis $w=w_0$, $m=m_0$, $f=f_0$, $b=b_0$ and $\rho=\rho_0$. In the remainder of this paper these variables have already been scaled so that $w_0=m_0=\rho_0=b_0=\rho_0=1$.

The continuity equation for the matrix density $\rho$ is
\begin{equation*}
\pd{\rho}{t}+\nabla\cdot (\rho \bv) = G_{\rho}(f,w,p),
\end{equation*}
where $G_{\rho}(f,w,\rho )$ is a growth and decay term of the ECM due to collagen secretion by fibroblasts
and degradation by matrix metalloproteinases (MMPs). The specific form of $G_{\rho}$ incorporates the
fact that collagen production and maturation require the availability of
oxygen \cite{Hutton:1967:CSR,Myllyla:1977:MPH,Gordillo:2003:RER,Sen:2009:WHE},
\begin{equation*}
G_{\rho}= \kineticrhon,
\end{equation*}
where $\rho_m$ is the maximum matrix volume fraction permitted in the partially healed region, $\rho_m>1$.

The partially healed tissue is modeled as a quasi-static upper convected Maxwell fluid with velocity $\bv$, deviatoric stress tensor given by
$\tau = \eta (\nabla \bv + \nabla \bv^{T})$, where $\eta$ is the shear viscosity, and pressure $P$.
The pressure $P$ is generally a function of the matrix density $\rho$, and is assumed to have the form
\beq\label{2.1}
P(\rho)=\left\{
\begin{aligned}
&\beta( \rho -1),  && \rho\geq 1 \\
&0, && \rho<1.
\end{aligned}\right.
\eeq
The total stress $\sigma = \tau -PI$   appears only in the boundary conditions.
By further assuming radially symmetric flow, i.e.,
$\bv = v(r,t)\mathbf{e}_r$, the continuity equation becomes
\beq\label{2.2}
\pd{\rho}{t}+\conv{\rho}= \kineticrhon, \quad R(t)<r<L,
\eeq
and the non-dimensionalized momentum equation for the ECM becomes (see \cite{Xue:2009:MMI}, supporting information)
\beq\label{2.3}
\frac{1}{r} \pd{}{r} \left( r \pd{v}{r} \right) - \frac{  v}{r^2}  = \pd{P(\rho)}{r}, \quad
R(t)<r<L.
\eeq

To simplify the analysis and simulations we wish to have a PDE system in which all variables are radially symmetric. In order to implement ischemic conditions
in radially symmetric form
we assume that small arcs of length $\delta$ are cut off from the healthy tissue at $r=L$ and that the distance between two adjacent $\delta$ arcs is $\varepsilon$. If $\delta, \; \varepsilon \to 0$ in such a way that $\eps\sim e^{-c/\delta}$ where $c$ is a positive constant, then, for any diffusion process with boundary conditions
\beaa
& \pd{u }{r} = 0 \;\; &\mbox{on the}\; \delta\mbox{-arcs}, \\
& u  = g   \;\;   &\mbox{on the remaining arcs},
\eeaa
the limiting ``homogenized'' boundary condition is \cite{Friedman:1995:EPB}
\beqn
(1-\gamma) (u-g) + \gamma  \pd{u}{r}  = 0 \on r=L
\eeqn
for some constant $\gamma\in[0,1]$ which depends only on $c$; $\gamma=0$ corresponds to healthy tissue (i.e., no excision of $\delta$-arcs) and $\gm$ near 1 corresponds to extreme ischemia.

The equations for the concentrations of oxygen, PDGF and VEGF are:
\bea
&& \frac{\partial w}{\partial t} + \conv{w} = \divp{D_w \pd{w}{r} }  \label{2.4}\\
&& \hspace*{2cm}+ \kineticwn, \nonumber \\
&&\frac{\partial p}{\partial t} + \conv{p} =   \divp{D_p \pd{p}{r}} +
\kineticpn, \label{2.5}  \\
&&\frac{\partial e}{\partial t} + \conv{e} =  \divp{D_e \pd{e}{r}  } +
\kineticen,
\label{2.6}
\eea

The equations for macrophages, fibroblasts, capillary tips and capillary sprouts include diffusion, generation and death of cells,
 and chemotactic migration of cells:
\bea
&&\hspace*{-0.4cm}\frac{\partial m}{\partial t} +\conv{m}  = \divp{ D_m \pd{m}{r}}  - \divp{
\bddtaxis{ \chi_m \rho  m H(1 - m/m_m) }{p} } \label{2.7}\\
&& \hspace*{5cm}  + \kineticmn,  \nonumber\\
&&\hspace*{-0.4cm} \frac{\partial f}{\partial t}+\conv{f} = \divp{D_f \pd{f}{r}} -\divp{
\bddtaxis{\chi_f \rho  f H(1-f/f_m) }{p}  } \label{2.8} \\
&& \hspace*{5cm} + \kineticfn,\nonumber
\eea
\bea
&& \frac{\partial n}{\partial t} + \conv{n}  =  \divp{ D_n \pd{n}{r}} - \divp{
\bddtaxis{\chi_n \rho n H(1-n/n_m)}{e} } \label{2.9}\\
&& \hspace*{4cm} + \kineticnn, \nonumber \\
&&\frac{\partial b}{\partial t} + \conv{b} = \divp{ D_b \pd{b}{r}} +
\divp{ \bddtaxis{A D_n  b }{n}} \label{2.10} \\
&&\hspace*{4cm} -
\divp{\bddtaxis{A \chi_n b \rho n H(1-n/n_m)}{e} }  \nonumber \\
&& \hspace*{4cm} +  \kineticbn.\nonumber
\eea
where the two terms with $A$ (in (\ref{2.10})) represent the fact that sprouts follow tips, and the oxygen-dependent functions $G$'s and $D$ are given by
\begin{align*}
&G_p(w) = \left\{
\begin{aligned}
& 3 w, && 0\leq w< 0.5  \\
&2- w, && 0.5 \leq w< 1  \\
& \frac{1}{3}w + \frac{2}{3},   && 1\leq w < 4  \\
& 2, &&  w \geq 4
\end{aligned}
\right., \qquad
G_e(w) = \left\{
\begin{aligned}
&2 w, && 0\leq w< 0.5, \\
&2- 2 w, && 0.5 \leq w< 1, \\
&\frac{1}{3}w - \frac{1}{3} , && 1 \leq w <4 , \\
& 1, && w \geq 4
\end{aligned}
\right.,
\end{align*}
\begin{align*}
&G_f(w) = \dfrac{(K_{wf} +1 ) w}{K_{wf}+w}, \quad G_b =
\dfrac{(K_{w\rho}+1) w}{K_{w\rho}+w}, \quad  D(w) = 1- H(5w-1)H(1 - w/3 ).
\end{align*}
Here $H$ is an approximated Heaviside function
\beqn
H(u) =\left\{
\begin{aligned}
&\frac{u^6}{10^{-6} + u^6},  && u\geq 0 \\
&0, && u<0  .
\end{aligned}
\right.
\eeqn

Note that in Equation \re{2.4} the supply of oxygen from the vasculature is reduced to $k_w b((1-\gm)w_b-w)$ due to the ischemic condition. The functions $G_p(w)$ and $G_e(w)$ are constructed to reflect the biological effect of oxygenation: moderate hypoxia and hyperoxia increase the production of PDGF and VEGF compared to normoxia. Equations \re{2.7} - \re{2.9} include chemotaxis flux terms that describe the chemotactic movement of macrophages, fibroblasts and capillary tips.
The two terms with $A$ in Equation \re{2.10} represent the fact that capillary sprouts are dragged along capillary tips.
Although the forms of the $G$ functions and $D$ function are suggested by biological experiments, our mathematical analysis will
not depend on the special form of these functions.

The free boundary $r=R(t)$ is moving with velocity $v$:
\begin{align}
\dot{R}(t) = v(R(t),t). \label{2.11}
\end{align}

The boundary conditions at $r=L$ are
\begin{align}
&v = 0, \label{2.12}\\
& (1-\gamma)(w-1) +   \gamma L  \pd{w}{r} =0, \label{2.13}
\end{align}
\begin{align}
&(1-\gamma)p +   \gamma L  \pd{p}{r} =0, \quad (1-\gamma)e +  \gamma L \pd{e}{r} =0, \label{2.14}\\
&(1-\gamma)m +   \gamma L \left(  \pd{m}{r} -  \dfrac{\chi_m }{D_m } \bddtaxis{\rho m H(1 - m/m_m)}{p}\right) =0, \label{2.15}\\
&(1-\gamma)(f-1) +  \gamma L  \left(  \pd{f}{r} -\dfrac{ \chi_f }{D_f } \bddtaxis{\rho f H(1-f/f_m)}{p} \right) =0,
\label{2.16}\\
&(1-\gamma)n +   \gamma L  \left(   \pd{n}{r} -   \dfrac{\chi_n }{D_n } \bddtaxis{\rho n H(1-n/n_m)}{e} \right) =0, \label{2.17} \\
&(1-\gamma)(b-1) +  \gamma L  \left( \pd{b}{r} + \bddtaxis{A D_n  b }{n} - \bddtaxis{A \chi_n b \rho n H(1-n/n_m)}{e} \right) = 0, \label{2.18}
\end{align}
and the boundary conditions at $r=R(t)$ are
\begin{align}
& \pd{v}{r} = P, \label{2.19}\\
& \pd{w}{r} = \pd{e}{r} =  \pd{n}{r}= \pd{b}{r} = 0, \label{2.20}   \\
& - \pd{p}{r} =  \frac{k_{pb} R(t)}{D_pR_0}, \label{2.21} \\
& - D_m \pd{m}{r} +  \chi_m \bddtaxis{\rho m H(1 - m/m_m)}{p} =  0,  \label{2.22}\\
& - D_f \pd{f}{r} +  \chi_f \bddtaxis{\rho f H(1-f/f_m)}{p}  =  0, \label{2.23}
\end{align}
Equation \re{2.21} represents the fact that secretion of platelets decreases with healing
(i.e., as $R(t)$ decreases).
The initial conditions for $R_0\leq r\leq L$ take the form
\beq\label{2.24}
\begin{aligned}
& R(0)= R_0, \;\; v = 0, \;\;  \rho   =   f  = 1, \;\;  w = 1,
\;\;  b  = g\left(\frac{r-R_0}{\epsilon_0}\right) ,\\
& e=m=n=0, \;\;
p=p_0(r),
\end{aligned}
\eeq
where
\begin{equation*}
g(z)=\left\{
\begin{aligned}
& 0, &&\qquad z\leq 0, \\
&\frac{8 }{3}z^2 , &&\qquad 0 < z \leq \frac{1}{4},\\
& \frac{4}{3}z - \frac{1}{6} , &&\qquad   \frac{1}{4}  \leq z < \frac{3}{4}  , \\
& 1 - \frac{8}{3 }(1-z)^2, &&\qquad  \frac{3}{4} \leq z \leq 1,\\
& 1, &&\qquad z > 1.
\end{aligned}\right.
\end{equation*}
and $p_0(r)$ has three continuous derivatives and satisfies the boundary conditions (\ref{2.14}) and (\ref{2.21}),
and
\beq\label{2.25}
\left\{
\begin{aligned}
& p_0 '(r) <0  &&\mbox{if}&& R_0 < r<R_0 + \eps_0, \\
& p_0 (r)=0 &&\mbox{if}&& R_0  + \eps_0 < r<L
\end{aligned}
\right.
\eeq
where $0<\eps_0<L-R_0$.

In a healthy tissue there is no net growth of ECM, i.e., $G_{\rho}(f,w,\rho)=0$ if $f=w=\rho=1$, which means that
\beq\label{2.26}
\lambda_{\rho} = \frac{k_{\rho}  }{1 +K_{w\rho}}\left( 1 -\frac{1}{\rho_m}\right).
\eeq
Similarly
\begin{equation}
k_w = \frac{ \ld{wf} + \ld{wm}  }{  w_b-1 },\label{2.27}
\end{equation}
\begin{equation}
k_f = \dfrac{\ld f}{1-1/f_m} . \label{2.28}
\end{equation}

\section{${R}(t)$ is monotonically decreasing}\label{3}

Set
\beq\label{3.1}
Q(t) =  \int_{R(t)}^L yP(y,t)dy,
\eeq
where $P(r,t) = P(\rho(r,t))$.

\begin{thm}
For any solution of (\ref{2.2})-(\ref{2.28}) there holds:
\bea
&& \dot{R}(t) \leq 0;  \; \;   \dot{R}(t) < 0  \; \mbox{ if and only if }
\;\; Q(t) > 0;   \label{3.2} \\
&&  R(0) e^{-\frac{2}{L^2}\int_0^t Q(\tau)d\tau} \leq R(t) \leq R(0) e^{-\frac{1}{L^2}\int_0^t Q(\tau)d\tau} \label{3.3}.
\eea
\end{thm}

\pf Equation (\ref{2.3}) can be rewritten as
\begin{equation*}
v_{rr} + \frac{ v_r }{r} -\frac{v }{r^2} =v_{rr} + \left(\frac{v}{r}\right)_r = P_r.
\end{equation*}
Integrating over $[R(t), r]$, we obtain
\begin{equation*}
v_r(r,t)-v_r(R(t))+\frac{ v(r,t)}{r} - \frac{v(R(t))}{R(t)} = P(r,t) - P(R(t),t).
\end{equation*}
From (\ref{2.11}) and (\ref{2.19}) we obtain
\begin{equation*}
v_r(r,t)+ \frac{ v(r,t)}{r} -  \frac{\dot{R}(t)}{R(t)} = P(r,t),
\end{equation*}
hence
\begin{equation}\label{3.4}
(rv)_r - r \frac{\dot{R}(t)}{R(t)} = rP(r,t).
\end{equation}
Integrating this equation over $[r, L]$ and using (\ref{2.12}), we obtain
\begin{equation}\label{3.5}
-r v(r,t) - \frac{L^2-r^2}{2}  \frac{\dot{R}(t)}{R(t)} = \int_r^L yP(y,t) dy.
\end{equation}
In particular, at $r=R(t)$,
\begin{equation*}
-R(t)\dot{R}(t)  - \frac{L^2-R(t)^2}{2}  \frac{\dot{R}(t)}{R(t)} = \int_{R(t)}^L yP(y,t) dy,
\end{equation*}
or
\begin{equation}
 \frac{\dot{R}(t)}{R(t)} = - \frac{2}{L^2+R(t)^2}Q(t).  \label{3.6}
\end{equation}
The assertion (\ref{3.2}) now follows immediately from (\ref{3.6}). From (\ref{3.6}) we also obtain
\beq\label{3.7}
-\frac{2}{L^2} Q(t)\leq \frac{\dot{R}(t)}{R(t)} \leq -\frac{1}{L^2} Q(t),
\eeq
from which we deduce the estimate (\ref{3.3}). \q

If we substitute $ \dot{R} /R $ from (\ref{3.6}) into (\ref{3.4}) we obtain, after dividing by $r$,
\beq\label{3.8}
\frac{(rv)_r}{r} = P(r,t) - \frac{2}{L^2+R(t)^2} Q(t);
\eeq
this equation will be needed in the remainder of this paper.
If we substitute $ \dot{R} /R $ from (\ref{3.6}) into (\ref{3.5}), and divide by $r$, we obtain an
expression for $v$,
\beqn
v(r,t) = \frac{1}{r}\left\{\frac{L^2-r^2}{L^2+R(t)^2} Q(t) - \int_r^L yP(y,t)dy\right\}
\eeqn
or
\beq\label{3.9}
v(r,t)=\frac{1}{r}\left\{ \frac{L^2-r^2}{L^2+R(t)^2} \int_{R(t)}^r yP(y,t)dy - \frac{r^2 +R(t)^2 }{L^2+R(t)^2}\int_r^L yP(y,t)dy \right\}
\eeq

\begin{cor}
Equation (\ref{2.3}) for $v$ together with the boundary conditions (\ref{2.12}), (\ref{2.19}) and the initial condition $v=0$ can be equivalently replaced by the formula (\ref{3.9}) .
\end{cor}

In the remainder of this paper we shall often work with the representation (\ref{3.9}) for $v$.

\section{ A priori estimates}\label{4}

In this section we assume that there exists a classical solution to (\ref{2.2}) -- (\ref{2.28}) for $0\leq t< T$, and derive \textit{a priori} estimates which depend on $T$, but remain uniformly bounded for any finite $T$. We set
\beqn
\Omega_{T}=\{ (r, \theta, t)\; |\; R(t)<r<L,\; 0 \leq \theta\leq 2\pi,\; 0<t \leq T\}
\eeqn
and introduce the following notation:

$C_{r,t}^{2+\alpha,1+\alpha/2}(\bar\Omega_T)$ is the space of functions $u(r,t)$ with $u$, $D^2_r u$, $D_t u$
uniformly H\"{o}lder continuous in $\bar\Omega_T$, with exponents
$\alpha$ in $r$ and $\alpha/2$ in t; the norm in this space is
defined by
\beqn
\| u \|_{C_{r,t}^{2+\alpha,1+\alpha/2}(\bar\Omega_T)} = \| u \|_{L^{\infty}(\bar\Omega_T)} + \| D^2_r u \|_{C_{r,t}^{\alpha,\alpha/2}(\bar\Omega_T)}+ \| D_t u \|_{C_{r,t}^{\alpha,\alpha/2}(\bar\Omega_T)}
\eeqn
where
\beqn
\| v \|_{C_{r,t}^{\alpha,\alpha/2}(\bar\Omega_T)} = \| v \|_{L^{\infty}(\Omega_T)} +
\sup_{(r,t), (r',t')) \in \bar\Omega_T} \frac{|v(r,t)-v(r',t')|}{|r-r'|^{\alpha}+|t-t'|^{\alpha/2}}.
\eeqn
Similarly we define the spaces
$C_{r,t}^{\alpha,\beta}(\bar\Omega_T)$,
$C^{1+\alpha}[0,T]$, etc.

\bigskip
In the remainder of this paper we shall use the following comparison principle \cite{FriedmanPPDE,LiebermanPDE}.

\begin{lma}\label{thm4.1}
Let $v_1$, $v_2$  satisfy
\beq\label{4.1}
\pd{v_1}{t} - D\Delta v_1 + g(x,t,v_1, \nabla v_1) \geq \pd{v_2}{t} - D\Delta v_2 + g(x,t,v_2,\nabla v_2) \iin \Omega_T.
\eeq
If
\beq\label{4.2}
\begin{aligned}
&  \mu_1\pd{}{\nu}(v_1-v_2) + \mu_2 (v_1-v_2) \geq 0 \on \partial \Omega_T\cap \{0<t<T\},\\
&  (v_1-v_2)|_{t=0} \geq 0 \iin \Omega_0
\end{aligned}
\eeq
where $\nu$ is the outward normal and $\mu_1$, $\mu_2$ are nonnegative
functions satisfying, at each point,
either $\mu_1>0$ or $\mu_1=0$,$\mu_2>0$, then $v_1\geq v_2$ in $\Omega_T$. Furthermore, if strict inequalities hold in both \re{4.1} and \re{4.2}, then 
$v_1>v_2$ in $\Omega_T$.
\end{lma}
\medskip

\begin{lma}\label{thm4.2}
For any solution of \re{2.2} -- \re{2.28},
\beq\label{4.3}
\mbox{the components } w,\; e,\; p,\; m,\; f,\; n,\; b,\; \mbox{ and } \rho \mbox{ are nonnegative functions.}
\eeq
\end{lma}

\vspace*{-0.4cm}
\pf For any small $\delta>0$, let us add $\delta$ on the right-hand side of each of the equations \re{2.4}-\re{2.10} and each of the boundary conditions \re{2.13}-\re{2.18}, \re{2.21}-\re{2.23}, replace $0$ by $-\delta$ in \re{2.20}, and increase the initial data of $b, e, m, n, p$ by $\delta$. We refer to this new system as the ``$\delta$-problem'' and to its solution as the ``$\delta$-solution''. By continuity, each component of the $\delta$-solution is strictly positive in $\Omega_{t_0}$ for some $t_0>0$. We claim that all the components are strictly positive in $\Omega_T$ for all $T>0$. Indeed, otherwise there is a smallest $T$ such that at least one component of the $\delta$-solution, denoted by $z$, vanishes at some point $(\bar{r}, T)$. We can then apply the second part of Lemma \ref{thm4.1}
with $v_1 = z$, $v_2 = 0$ to conclude that $z(\bar{r},T)>0$, which is a contradiction.

The local existence and uniqueness proof given in Sections 4-6 is valid also for the $\delta$-problem. The estimates derived there are uniform in $\delta$ so that, as $\delta\to 0$, the $\delta$-solution converges to the original solution. Hence each component of the original solution is non-negative in a small time interval, say $0<t<t_{\ast}$. We can now repeat the process for $t>t_{\ast}$, and conclude, step-by-step that each component of the solution is non-negative in $\Omega_T$ for any $T>0$. \q

\begin{lma}\label{thm4.3}
If initially $\rho(r,0)<\rho_m$ for $R(0)\leq r\leq L$, then, 
\begin{align}
\rho<\rho_m  \qquad \mbox{in}\;\;   \Omega_T. \label{4.4}
\end{align}
\end{lma}
\pf
If the assertion (\ref{4.4}) is not true, then there exists a
$t^{\ast}>0$ such that
$\rho(r,t)<\rho_m$ in $\Omega_{t^{\ast}}$,
and $\rho(r^{\ast}, t^{\ast}) = \rho_m$ for some $R(t^{\ast}) \leq r^{\ast} \leq L$.
Then, along the characteristic curve with velocity $v$, through $(r^{\ast},t^{\ast})$,
\beq\label{4.5}
\left.\frac{D\rho}{Dt}\right|_{(r^{\ast}, t^{\ast})}  \geq 0,
\eeq
where $D/Dt = \partial /\partial t + v(\partial/\partial r) $.
On the other hand, from (\ref{2.2}) and (\ref{3.8}) we get,
\begin{align*}
\left.\frac{D\rho}{Dt}\right|_{(r^{\ast}, t^{\ast})} = -\ld{\rho}\rho(r^{\ast}, t^{\ast}) - \left(P(r^{\ast}, t^{\ast}) - \frac{2}{L^2+R^2}Q(t^{\ast})\right)\rho(r^{\ast}, t^{\ast}).
\end{align*}
Since $Q(t^{\ast})\leq \frac{L^2-R^2}{2} \max_r P(r,t^{\ast}) = \frac{L^2-R^2}{2} P(r^{\ast}, t^{\ast})$, we obtain
\begin{align*}
\left.\frac{D\rho}{Dt}\right|_{(r^{\ast}, t^{\ast})} = -\ld{\rho}\rho(r^{\ast}, t^{\ast}) - \frac{2 R^2}{L^2+R^2}P(r^{\ast}, t^{\ast}) \rho(r^{\ast}, t^{\ast}) <0.
\end{align*}
which is a contradiction to (\ref{4.5}). \q

Recall that we have assumed $\rho_m > 1$. 
\begin{lma} \label{thm4.4}
There holds:
\beq
\frac{\left| v(r,t)  \right|}{r} \leq \beta(\rho_m -1), \quad |v_r(r,t)| \leq 2\beta(\rho_m-1), \iin  \Omega_T.
\eeq
\end{lma}
\pf
From Lemma  \ref{thm4.3} we obtain
\beaa
\int_{R(t)}^r y P(y,t)dy \leq \beta(\rho_m - 1) \frac{r^2-R(t)^2}{2}, \\
 \int_r^L y P(y,t)dy \leq \beta(\rho_m - 1) \frac{L^2 - r^2}{2}.
\eeaa
Using these estimates in (\ref{3.9}) we get
\beaa
\frac{|v(r,t)|}{r} \leq  \beta(\rho_m - 1) \frac{L^2 - r^2}{L^2 + R(t)^2} \leq \beta(\rho_m - 1) .
\eeaa
Substituting this inequality into (\ref{3.8}) and estimating $P$ and $Q$ by Lemma \ref{thm4.3}, we also obtain
\beqn
\left| v_r(r,t) \right| \leq 2 \beta(\rho_m - 1). \q
\eeqn

\begin{lma}\label{thm4.5}
Setting \begin{align*}
N = \max\left\{ \frac{k_{nb}}{\ld{nb}}, \, \frac{k_n + \beta [\rho_m - 1]}{\ld{nn}},\, n_m \right\},
\end{align*}
there holds:
\begin{align}
0\leq n(r,t)\leq N  \qquad \mbox{in}\;\;   \Omega_T.  \label{4.7}
\end{align}
\end{lma}
\pf
We write Equation (\ref{2.9}) for $n$ in the form
\begin{align*}
\mathscript{L}[n] = \mathscript{L}_0[n] + \mathscript{F}[n] = 0,
\end{align*}
where
\begin{align*}
\sL_0[\phi]=\pd{\phi}{t} - \frac{1}{r} \pd{}{r} \left(rD_n\pd{\phi}{r}\right) + v \phi_r +
\frac{1}{r}\pd{}{r}
\left( r\, \bddtaxis{\chi_n\rho\phi H(1-\phi/n_m)} {e} \right),
\end{align*}
and
\beqn
\sF[\phi] = b\left(\ld{nb}\phi - k_{nb} \frac{e}{1+e}\right) + \left(\ld{nn}\phi + \frac{(rv)_r}{r} -k_n \frac{e}{1+e}\right)\phi.
\eeqn
By (\ref{3.8}) and Lemma \ref{thm4.3}
\begin{align}
\frac{1}{r}(rv)_r  \geq -\beta [\rho_m  -1 ],\nonumber
\end{align}
so that, by definition of $N$,
\beqn
\ld{nn} N + \frac{(rv)_r}{r} - k_{n} \frac{e}{1+e} > \ld{nn} N -\beta [\rho_m  -1 ]- k_{n} >0,
\eeqn
and
\beqn
\ld{nb}N - k_{nb}\frac{e}{1+e} > \ld{nb}N - k_{nb} \geq 0.
\eeqn
Since, by (\ref{4.3}), $b\geq 0$, we conclude that $\sF[N]\geq 0$ and hence $N$ is a supersolution, i.e., $\mathscript{L}(N)\geq 0$. Using also the boundary conditions (\ref{2.17}) and (\ref{2.20}) we deduce, by the comparison lemma, that $n(r,t)\leq N.$ \q

\begin{lma}\label{4.6}
For any $T>0$, there exists a constant $C_T$ such that
\beq\label{4.8}
0\leq b(r,t)\leq C_T   \qquad \mbox{in}\;\;   \Omega_T     .
\eeq
\end{lma}
\pf
By the comparison principle,
\beqn
0\leq b(r,t)\leq b_1(r,t)
\eeqn
where $b_1(r,t)$ is a solution of the same equation as $b(r,t)$ but without the quadratic term $-k_b G_b(w) b^2$ and with the same boundary and initial conditions as for $b(r,t)$. We can write the equation for $b_1$ in the form
\beq\label{4.9}
r\pd{b_1}{t}-\pd{}{r}\left(rD_b\pd{b_1}{r}\right) +a_1(r,t)b_1(r,t) + a_2(r,t)\pd{b_1}{r}(r,t) + \pd{}{r}\big(a_3(r,t)b_1(r,t)\big)=a_4(r,t),
\eeq
where, by using (\ref{4.7}), we find that $a_1$, $a_2$, $a_3$, $a_4$ are all uniformly bounded.
From the Nash-Moser estimate \cite{LiebermanPDE} we deduce that, for any $0<t_1 \leq T$,
\beq\label{4.10}
\|b_1\|_{C^{\alpha, \alpha/2}(\Omega_{t_1})} \leq C_T + C_T \|b_1\|_{L^{\infty}(\Omega_{t_1})} ,
\eeq
and by interpolation,
\begin{eqnarray*}
\|b_1\|_{L^{\infty}(\Omega_{t_1})} &\leq& \|b_1(\cdot,0)\|_{L^{\infty} } + t_1^{\alpha/2}(1+\sup_{0\leq \tau\leq t_1}{|\dot{R}(t)|^{\alpha/2}}) \|b_1\|_{C^{\alpha, \alpha/2} }\\
&\leq & \|b_1(\cdot,0 )\|_{L^{\infty}}  + C^{\ast} t_1^{\alpha/2}\big(C_T + C_T \|b_1\|_{L^{\infty}(\Omega_{t_1})}\big)\\
&\leq &  C^{\ast} C_T t_1^{\alpha/2} \|b_1\|_{L^{\infty}(\Omega_{t_1})} + C.
\end{eqnarray*}
Choosing $t_1$ such that $$C^{\ast} C_T t_1^{\alpha/2} = \frac{1}{2}, $$ we obtain the estimate
$$ \|b_1\|_{L^{\infty}(\Omega_{t_1})} \leq C. $$
Repeating this procedure step-by-step, the assertion (\ref{4.8}) follows. \q

The above proof can be applied successively to $m$, $f$,  $p$, $e$ and $w$ to establish the following estimates.

\begin{lma}\label{thm4.7}
For any $T>0$, there exists a positive constant $C_T$ such that $\iin \Omega_{T}$,
\bea
&0\leq m(r,t)\leq C_T, \;\; 0\leq f(r,t)\leq C_T, \;\; 0\leq p(r,t)\leq C_T, \nonumber \\
&0\leq e(r,t)\leq C_T, \;\; 0\leq w(r,t)\leq C_T  . \label{4.11}
\eea
\end{lma}

Since $b$ is bounded (by $C_T$) in $\Omega_{T}$, we can write the equation (\ref{2.10}) for $b$ in the same form as Equation (\ref{4.9}) for $b_1$ and thus derive, by the Nash-Moser estimate, a H\"{o}lder bound
\beqn
\|b\|_{C^{\alpha, \alpha/2}(\bar\Omega_T)} \leq C_T.
\eeqn
The same bound can be derived for the components $n$, $m$, $f$, $p$, $e$ and then also for $w$. Hence, we obtain

\begin{lma}\label{thm4.8}
For any $T>0$ there exists a positive constant $C_T$ such that
\beq
\|w, p, e, m, f, n, b\|_{C^{\alpha, \alpha/2}(\bar\Omega_T)}\leq C_T.
\eeq
\end{lma}

\bigskip
Rewriting Equation (\ref{2.2}) in the form
\beq\label{4.13}
\rho_t + v \rho_r = \kineticrhon - \frac{(rv)_r}{r}\rho \equiv \cF(r,s),
\eeq
we proceed to establish a H\"older estimate for the function $\rho$.

\begin{lma}\label{thm4.9}
For any $T>0$ there exists a constant $C_T$ such that
\beq\label{4.14}
\| \rho \|_{C_{r,t}^{\alpha,\alpha}(\bar\Omega_T)} \leq C_T.
\eeq
\end{lma}
\pf We introduce the characteristic curves $X$, for (\ref{4.13}), by
\beqn\left\{
\begin{aligned}
&\dd{X_r(r,t,s)}{s} = v_r(X(r,t,s),s)X_r(r,t,s), \:\: \forall s\in [0,t]\\
&X_r(r,t,t) = 1.
\end{aligned}\right.
\eeqn
Using Lemma \ref{thm4.4} we find that
\beqn
|X_r(r,t,s)| \leq e^{2 \beta(\rho_m - 1)(t-s)}.
\eeqn

Let $J(r,t,s) = \rho(X(r,t,s),s)$, so that
\beqn\left\{
\begin{aligned}
&\dd{J(r,t,s)}{s} = \cF(X(r,t,s),s),  \\
&J(r,t,t) = \rho(r,t).
\end{aligned}\right.
\eeqn
Then
\beaa
&& \hspace*{-1.5cm} \frac{|\rho(r_1,t) - \rho(r_2,t)|}{|r_1 - r_2|^{\alpha}} \\
&\leq& \frac{1}{|r_1 - r_2|^{\alpha}} \left| \int_0^t \cF(X(r_1,t,s),s )- \cF(X(r_2,t,s),s ) ds\right| \\
&& +  \frac{ |\rho(X(r_1,t,0), 0) - \rho(X(r_2, t, 0), 0)| }{|r_1 - r_2|^{\alpha}}. \\
\eeaa
By the initial condition $\rho(r,0) \equiv 1$ the last term vanishes, and
\beaa
&& \hspace*{-.5cm} \frac{1}{|r_1 - r_2|^{\alpha}} \left| \int_0^t \cF(X(r_1,t,s),s ) - \cF(X(r_2,t,s),s ) ds\right|\\
&\leq& \left| \int_0^t \frac{\cF(X(r_1,t,s),s ) - \cF(X(r_2,t,s),s)}{|X(r_1,t,s)-X(r_2,t,s))|^{\alpha}} \cdot
\left(\frac{|X(r_1,t,s)-X(r_2,t,s)|}{|r_1 - r_2|}\right)^{\alpha} ds\right| \\
&\leq& (e^{2  \beta(\rho_m - 1)(t-s)})^{\alpha} \int_0^t \frac{|\cF(X(r_1,t,s),s ) - \cF(X(r_2,t,s),s)|}{|X(r_1,t,s)-X(r_2,t,s))|^{\alpha}} ds \\
&\leq& C_T \int_0^t   [\rho(\cdot,s)]_{C_r^{\alpha}} + [w(\cdot,s)]_{C_r^{\alpha}}  + [f(\cdot,s)]_{C_r^{\alpha}}  ds.
\eeaa
Hence
\beqn
\frac{|\rho(r_1,t) - \rho(r_2,t)|}{|r_1 - r_2|^{\alpha}} \leq C_T + C_T \int_0^t  [\rho(\cdot,s)]_{C_r^{\alpha}} ds.
\eeqn
Taking supremum over  $r_1, \; r_2 \in [R(t),\, L]$, $r_1\neq r_2$, we obtain
\beqn
 [\rho(\cdot,t)]_{C_r^{\alpha}} \leq C_T + C_T \int_0^t  [\rho(\cdot,s)]_{C_r^{\alpha}} ds,
\eeqn
and by Gronwall's inequality,
\beq\label{4.15}
[\rho(\cdot,t)]_{C_r^{\alpha}} \leq C_T.
\eeq

Next, taking $t_2>t_1>0$,  we can write
\beaa
 \rho(r,t_2) -\rho(r, t_1)
 = \int_{t_1}^{t_2} \cF(X(r,t_2,s),s)ds + \rho(X(r,t_2,t_1),t_1) - \rho(r,t_1),
\eeaa
so that
\beqn
 \rho(r,t_2) -\rho(r, t_1) \leq C |t_2-t_1| + [\rho(\cdot,t_1)]_{C_r^{\alpha}} |X(r,t_2,t_1)-r|^{\alpha}.
\eeqn
Since
\beqn
|X(r,t_2,t_1)-r| = |X(r,t_2,t_1)- X(r,t_2,t_2) | \leq \left\| \dd{X}{s} \right\|_{L^{\infty}} |t_2-t_1|,
\eeqn
we obtain
\beqn
|\rho(r,t_2) -\rho(r, t_1)|\leq C_T|t_2-t_1|^{\alpha}.
\eeqn
Combining this inequality with (\ref{4.15}), the assertion (\ref{4.14}) follows. \q

\begin{lma}\label{thm4.10}
For any $T>0$ there exists a constant $C_T$ such that
\beq
\|v\|_{C_{r,t}^{\alpha,\alpha}(\bar\Omega_T)} + \|v_r\|_{C_{r,t}^{\alpha,\alpha}(\bar\Omega_T)}\leq C_T.\label{4.16}
\eeq
\end{lma}
\pf The proof follows from the representations of $v(r,t)$ and $v_r(r,t)$ in (\ref{3.9}) and (\ref{3.8}) by using Lemma \ref{thm4.9} and the boundedness of $\dot{R}$ (from (\ref{3.3})). \q

\begin{lma}\label{thm4.11}
For any $T>0$ there exists a constant $C_T$ such that
\beq\label{4.17}
\|R\|_{C^{1+\alpha}([0,T])} \leq C_T. \eeq
\end{lma}
\pf This follows from (\ref{2.11}) and Lemma \ref{thm4.10}. \q

\begin{lma}\label{4.12}
For any $T>0$ there exists a constant $C_T$ such that

(i)
\beaa
&&\|p\|_{C_{r,t}^{2+\alpha,1+\alpha/2}( \bar{\Omega}_T)} \leq C_T, \\
&&\|e\|_{C_{r,t}^{2+\alpha,1+\alpha/2}( \bar{\Omega}_T)} \leq C_T, \\
&&\|w\|_{C_{r,t}^{2+\alpha,1+\alpha/2}( \bar{\Omega}_T)} \leq C_T; \\
\eeaa

(ii)
\beaa
&&\|m \|_{C_{r,t}^{2+\alpha,1+\alpha/2}( \bar{\Omega}_T)} \leq C_T, \\
&&\|f\|_{C_{r,t}^{2+\alpha,1+\alpha/2}( \bar{\Omega}_T)} \leq C_T, \\
&&\|  n \|_{C_{r,t}^{2+\alpha,1+\alpha/2}( \bar{\Omega}_T)} \leq C_T, \\
&&\|  b\|_{C_{r,t}^{2+\alpha,1+\alpha/2}( \bar{\Omega}_T)} \leq C_T;
\eeaa

(iii)
\beaa
&& \| \rho \|_{C_{r,t}^{2+\alpha,1+\alpha/2}( \bar{\Omega}_T)} \leq C_T, \\
&& \| v \|_{C_{r,t}^{2+\alpha,1+\alpha/2}( \bar{\Omega}_T)} \leq
C_T.
\eeaa
\end{lma}
\pf Indeed, (i) follows from Lemmas \ref{thm4.8} -- \ref{thm4.11} and the parabolic Schauder estimates \cite{FriedmanPPDE, LiebermanPDE}. The assertion (ii) follows by the Schauder estimates and (i). To prove (iii) we first formally differentiate (\ref{4.13}) in $r$ and apply the proof of Lemma \ref{thm4.9}, making use of Lemma \ref{thm4.10} and (ii). We thus obtain the bound
\beq\label{4.18}
\| \rho_r \|_{C_{r,t}^{\alpha,\alpha/2}( \bar{\Omega}_T)} \leq C_T.
\eeq

In order to rigorously prove \re{4.18}, we consider the solution $\tilde{\rho}_r$ of the differentiated equation (\ref{4.13}) and derive the estimate (\ref{4.18}). By integration of the equation of $\tilde{\rho}_r$ with respect to $r$, one can verify that $\int^r \tilde{\rho}_r dr$ coincides with $\rho$; hence $\partial \rho /\partial r = \tilde{\rho}_r$ and (\ref{4.18}) follows.

Differentiating (\ref{3.8}) in $r$ and using (\ref{4.18}) we deduce that
\beqn
\| v_{rr} \|_{C_{r,t}^{\alpha,\alpha/2}( \bar{\Omega}_T)} \leq C_T.
\eeqn
and this allows us to differentiate the equation for $\rho_r$ once more in $r$. Proceeding as before it is then easy to complete the proof of (iii).

\section{Transformation to a fixed domain}\label{5}

In order to prove existence and uniqueness of a solution of (\ref{2.2}) -- (\ref{2.28}) for a small time interval $0<t<T$, it is convenient to transform the system with the free boundary $r=R(t)$ into a system with a fixed boundary, using the mapping
\begin{equation}
\xi = \frac{r-R(t)}{L-R(t)}, \quad \big(r=(1-\xi)R(t) + \xi L\big).    \label{5.1}
\end{equation}
In the new system $\xi$ varies in the interval $0<\xi<1$, and for any function $u(r,t)=\tilde{u}(\xi,t)$,
\bea
&& \frac{\partial u}{\partial r}  =
\frac{1}{L-R(t)} \frac{\partial \tilde{u}}{\partial \xi} , \label{5.2} \\
&&  \frac{\partial }{\partial r} \left(r  \frac{\partial u}{\partial r} \right)
=\frac{1}{(L-R(t))^2} \frac{\partial }{\partial \xi}\left( r(\xi) \frac{\partial \tilde{u}}{\partial \xi }\right),\label{5.3}
\eea
and
\beaa
&&  \dfrac{\partial u}{\partial t} = \dfrac{\partial \tilde{u}}{\partial t}+\dfrac{\partial \tilde{u}}{\partial \xi}\dfrac{\partial \xi}{\partial t} = \frac{\partial \tilde{u}}{\partial t} + \frac{ \dot{R}(t)}{L-R(t)} \big(\xi - 1\big)  \frac{\partial \tilde{u}}{\partial \xi},  \\
\eeaa
\beaa
&& \big(\xi-1\big)\frac{\partial \tilde{u}}{\partial \xi} = \frac{1}{r}\pd{}{\xi}\bigg( r \big(\xi-1\big) \tilde{u} \bigg) + \left(\frac{(1- \xi )(L-R(t))}{r} - 1 \right) \tilde{u}.
\eeaa
Using these formulas we compute
\begin{equation*}
 \dfrac{\partial u}{\partial t}  + \conv{u} =  \dfrac{\partial \tilde{u}}{\partial t}  + B,
\end{equation*}
where
\begin{equation*}
\begin{aligned}
B  &=&& \frac{ \dot{R}(t)}{L-R(t)} \big(\xi - 1\big)  \frac{\partial \tilde{u}}{\partial \xi} + \frac{1}{(L-R(t))r} \pd{}{\xi}\big(r\tilde{u}v\big), \\
 &=&& \frac{ \dot{R}(t)}{L-R(t)} \frac{1}{r}\pd{}{\xi}\bigg( r \big(\xi-1\big) \tilde{u} \bigg)  + \frac{1}{(L-R(t))r} \pd{}{\xi}\big(r\tilde{u}v\big) + K \tilde{u},
\end{aligned}
\end{equation*}
or,
\beaa
B= \frac{1}{L-R(t)}\left[ \frac{1}{r} \pd{}{\xi}\Big( r \tilde{u} \big(\dot{R}(t) (\xi - 1) + v \big)  \Big)\right] + K \tilde{u}.
\eeaa
where
\begin{equation}
K=K(\xi) = \dfrac{\dot{R}(t)}{L-R(t)} \left(\dfrac{(1-\xi)(L-R(t))}{r} -1 \right), \label{5.4}
\end{equation}
Hence
\begin{equation}
\dfrac{\partial u}{\partial t} + \conv{u} = \pd{\tilde{u}}{t} + \frac{1}{L-R(t)}\left[ \frac{1}{r} \pd{}{\xi}\Big( r\tilde{u}  \big(\dot{R}(t) (\xi - 1) + v \big)  \Big)\right] + K \tilde{u}. \label{5.5}
\end{equation}
Using (\ref{5.2}), (\ref{5.3}) and (\ref{5.5}), we can transform the PDEs in Section \ref{2} into the following system of equations, where we have, for simplicity, dropped the tilda ``$\sim$'' from all the variables:
\begin{align}
& \pd{\rho}{t} +   \frac{1}{r(\xi)} \pd{}{\xi}\Big( r(\xi) \rho  M  \Big)    = \kineticrhon - K \rho, \hspace*{1cm} \label{5.6}\\
& \frac{1}{(L-R(t))^2}\frac{1}{r(\xi)}\pd{}{\xi}\left(r(\xi)\pd{v}{\xi}\right) -\frac{v}{r^2(\xi)} = \frac{1}{L-R(t)}\pd{P}{\xi}, \label{5.7} \\
& \pd{w}{t} +    \frac{1}{r(\xi)} \pd{}{\xi}\Big( r(\xi) w M  \Big)   =   \frac{1}{r(\xi)} \pd{}{\xi}\left(r(\xi)    \diffct{w}\pd{w}{\xi} \right)  \\ & \hspace*{1cm} + \kineticwn - K w,  \nonumber \\
& \pd{p}{t} +  \frac{1}{r(\xi)} \pd{}{\xi}\Big( r(\xi) pM  \Big)  = \frac{1}{r(\xi)} \pd{}{\xi}\left(r(\xi)\diffct{p} \pd{p}{\xi} \right)  \\ & \hspace*{3cm}  + \kineticpn - K p, \nonumber \\
& \pd{e}{t} + \frac{1}{r(\xi)} \pd{}{\xi}\Big( r(\xi) e M \Big)  =  \frac{1}{r(\xi)} \pd{}{\xi}\left(r(\xi)\diffct{e}\pd{e}{\xi} \right) \\ & \hspace*{3cm} + \kineticen  - K e, \nonumber\\
\end{align}
\begin{align}
&  \pd{m}{t} +\frac{1}{r(\xi)} \pd{}{\xi}\Big( r(\xi) m M \Big)  \\ & \hspace*{1cm} = \frac{1}{r(\xi)} \pd{}{\xi}\Big(r(\xi)  \diffct{m} \pd{m}{\xi}\Big)   -
\frac{1}{r(\xi)} \pd{}{\xi}\left(r(\xi)\bddtaxist{ \chit{m} \rho m H(1-m/m_m)}{p} \right)  \nonumber \\ & \hspace*{3cm} +  \kineticmn - K m , \nonumber\\
&  \pd{f}{t} +  \frac{1}{r(\xi)} \pd{}{\xi}\Big( r(\xi) f M \Big) \\ & \hspace*{1cm} = \frac{1}{r(\xi)} \pd{}{\xi}\Big(r(\xi)  \diffct{f} \pd{f}{\xi}\Big) -
 \frac{1}{r(\xi)} \pd{}{\xi}\left(r(\xi) \bddtaxist{\chit{f} \rho f H(1-f/f_m)}{p}  \right) \nonumber \\ & \hspace*{3cm} +  \kineticfn - K f \nonumber, \\
& \pd{n}{t} +  \frac{1}{r(\xi)} \pd{}{\xi}\Big( r(\xi) n M \Big) \\ & \hspace*{1cm} =  \frac{1}{r(\xi)} \pd{}{\xi}\Big(r(\xi)  \diffct{n} \pd{n}{\xi}\Big) - \frac{1}{r(\xi)} \pd{}{\xi}\left(r(\xi)  \bddtaxist{\chit{f} \rho n H(1-n/n_m)}{e}   \right)  \nonumber \\ & \hspace*{3cm} + \kineticnn - K n , \nonumber  \\
& \label{5.15}\pd{b}{t} +  \frac{1}{r(\xi)} \pd{}{\xi}\Big( r(\xi) b M \Big)  \\ &  \hspace*{1 cm} =  \frac{1}{r(\xi)} \pd{}{\xi}\Big(r(\xi)  \diffct{b} \pd{b}{\xi}\Big)  + \frac{1}{r(\xi)} \pd{}{\xi}\left(r(\xi)  \bddtaxist{A \diffct{n} b}{n}\right)\nonumber \\ &  \hspace*{2 cm}  - \frac{1}{r(\xi)} \pd{}{\xi}\left(r(\xi)  \bddtaxist{A \chit{n} b \rho n H(1-n/n_m)}{e}   \right)\nonumber \\ &  \hspace*{3 cm}  +\kineticbn - K b , \nonumber
\end{align}
where
\begin{align*}
&M = \dfrac{\dot{R}(t) (\xi - 1) + v  }{L-R(t)} , \qquad k_{sg}(t)=\frac{k_{sg}}{(L-R(t))^2},\\
&\diffct{u}=\frac{D_u}{(L-R(t))^2}, \quad \mbox{for } u = w, p, e, m, f, n, b, \\
&\chit{u}=\frac{\chi_u}{(L-R(t))^2}, \quad \mbox{for } u =  m, f, n, b.\\
\end{align*}
The free boundary condition remains as before, namely,
\beq
\dot{R}(t) = v(R(t),t). \label{5.16}
\eeq

The boundary conditions at the fixed boundary $\xi = 1$ are
\begin{align}
& v = 0, \label{5.17}\\
& (1-\gamma)(w-1) + \frac{\gamma  L }{L-R(t)} \pd{w}{\xi} =0   , \\
& (1-\gamma)p +  \frac{\gamma L }{L-R(t)} \pd{p}{\xi} =0 ,\\
& (1-\gamma)e + \frac{\gamma L }{L-R(t)} \pd{e}{\xi} =0, \\
& (1-\gamma)m + \frac{\gamma L }{L-R(t)} \left( \pd{m}{\xi} - \frac{\chi_m}{D_m} \bddtaxist{\rho m H(1 - m/m_m)}{p}\right) =0, \\
&  (1-\gamma)(f-1) + \frac{\gamma L }{L-R(t)} \left( \pd{f}{\xi} - \frac{\chi_f}{D_f} \bddtaxist{\rho f H(1-f/f_m)}{p} \right) =0  ,\\
& (1-\gamma)n + \frac{\gamma L }{L-R(t)} \left( \pd{n}{\xi} - \frac{\chi_n}{D_n}\bddtaxist{\rho n H(1-n/n_m) }{n} \right) =0  , \\
& (1-\gamma)(b - 1) \\
& \hspace*{0.3cm}+ \frac{\gamma L}{L-R(t)}\left( \pd{b}{\xi}  + \bddtaxist{A D_n  b }{n} - \bddtaxist{A \chi_n b \rho n H(1-n/n_m)}{e} \right)  = 0, \nonumber
\end{align}
and at the free boundary $\xi = 0$ they are
\begin{align}
&  \pd{v}{\xi} = \big(L- R(t)\big)P , \\
&  \pd{w}{\xi} = \pd{e}{\xi} = \pd{n}{\xi}= \pd{b}{\xi}= 0, \\
&  \pd{p}{\xi} =  - \frac{k_{pb} R}{D_p R_0} \big(L- R(t)\big),\\
& - D_m \pd{m}{\xi} +  \chi_m \bddtaxist{\rho m H(1-m/m_m) }{p} =  0, \\
& - D_f \pd{f}{\xi} +  \chi_f \bddtaxist{\rho f H(1-f/f_m)}{p} =  0.
\end{align}
The initial conditions take the form
\beq\label{5.30}
\begin{aligned}
& R(0)=R_0, \;v = 0, \; \rho  = f = 1 , \;  w = 1, \;  b = g\left(\frac{\xi(L-R_0)}{\eps_0}\right) ,\\
& e=m=n=0, \; p(\xi,0)= p_0\big((1-\xi)R_0 + \xi L\big).
\end{aligned}
\eeq

\section{Existence and Uniqueness}\label{6}
In this section we prove the following theorem.
\begin{thm}\label{6.1}
There exists a unique solution of (\ref{2.2}) -- (\ref{2.28})  for $0\leq t< \infty$ such that, for each $T>0$, the estimates of Lemma \ref{4.12} hold.
\end{thm}
\pf
We first prove existence and uniqueness for a small time interval $0\leq t\leq \tau$. For this proof it will be convenient to transform the system (\ref{2.2}) -- (\ref{2.24}) into the system (\ref{5.6}) -- (\ref{5.30}) with a fixed boundary.
Set
\beqn
G = \{0\leq \xi\leq 1\}, \quad G_T = \{(\xi,t); \xi\in G, 0\leq t\leq T\} \quad \mbox{for any} \;\;\; T>0,
\eeqn
and introduce the Banach space
\beaa
Y = \{(R(t), \rho(\xi, t)); R(0)= R_0, \rho(\xi,0)=1 \;\;\mbox{with norm} \\
\|(R, \rho) \|_Y = \|R\|_{C^{1+\alpha/2}[0,\tau]} + \|(\rho, \rho_{\xi})\|_{C^{\alpha,\alpha/2}(\bar G_{\tau})} \}
\eeaa
and the ball
\beqn
Y_B = \{(R,\rho)\in Y; \| (R, \rho)\|_Y \leq B\}
\eeqn
for any $B> 1+R_0$.

For any $(R,\rho) \in Y_B$ we wish to solve the system (\ref{5.7}) -- (\ref{5.15}) with the corresponding boundary and initial conditions from (\ref{5.17}) -- (\ref{5.30}). Denoting this solution by $u=(w,p,e,m,f,n,b,v)$ we shall then define $(\tilde{R},\tilde{\rho})$ by
\beq\label{thm6.1}
\frac{d}{dt}\tdR(t) = v(R(t),t), \quad \tdR(0) = R_0,
\eeq
\beq\label{6.2}
\pd{\tdrho}{t} + \frac{1}{r(\xi)}\pd{}{\xi}\Big(r(\xi)\tdrho\tdM\Big) = \frac{k_{\rho} w }{w+K_{w\rho}} f (1-\frac{\tdrho}{\rho_m}) - \lambda_{\rho} \tdrho -\tdK\tdrho, \quad \tdrho(\xi,0)=1,
\eeq
where
\beqn
\tdM = \dfrac{(d\tdR/dt) (\xi - 1) + v  }{L-\tdR(t)}, \quad \tdK=\dfrac{d\tdR/dt}{L-\tdR(t)} \left(\dfrac{(1-\xi)(L-\tdR(t))}{r} -1 \right),
\eeqn
and set
\beqn
(\tdR, \tdrho) = W(R,\rho)
\eeqn
We aim to prove that the mapping $W$ is a contraction mapping, and thus has a unique fixed point.

As in \cite{Friedman:2005:AMM} one can prove, by a fixed point argument, that there exists a unique solution $u$ for $0\leq t\leq \tau$, for $\tau$ small, and that
\beq\label{6.3}
\|u\|_{C_{\xi,t}^{2+\alpha, 1+\alpha/2}(\bar G_{\tau})}\leq C, \quad u=(w,p,e,m,f,n,b,v).
\eeq
The estimate (\ref{6.3}) can also be established by the argument used in the proof of Lemma \ref{4.12}. From (\ref{6.1}) and (\ref{6.3}) we get
\beq\label{6.4}
\|\frac{d}{dt}\tdR\|_{C^{2+\alpha}[0,\tau]} \leq C,
\eeq
so that
\beqn
\|(\tdM, \tdK)\|_{C_{\xi,t}^{2+\alpha, 1+\alpha/2}(\bar G_{\tau})} \leq C.
\eeqn
We next consider (\ref{6.2}), and use the same arguments as in the proofs of Lemma \ref{thm4.9} and \ref{4.12} (iii), to derive the estimate
\beq\label{6.5}
\|\tdrho\|_{C_{\xi,t}^{2+\alpha, 1+\alpha/2}(\bar G_{\tau})} \leq C.
\eeq
From (\ref{6.4}), (\ref{6.5}) we deduce that
\beq\label{6.6}\left\{
\begin{aligned}
&\|\tdR\|_{C^{1+\alpha}[0,\tau]}\leq R_0+C\tau, \\
&\|(\tdrho, \tdrho_{\xi})\|_{C_{\xi,t}^{\alpha, \alpha/2}(\bar G_{\tau})}\leq 1+C\tau^{1/2}.
\end{aligned}\right.
\eeq
Hence if $\tau$ is sufficiently small then $W$ maps $Y_B$ into itself.

We next prove that $W$ is a contraction in $Y_B$. Let $(R_1, \rho_1)$ and  $R_2, \rho_2$ be any elements in $Y_B$ and denote the corresponding solution by $u_i = (w_i, p_i, e_i, m_i, f_i, n_i, b_i, v_i)$ for $i=1, 2$. Set
\beqn
(\tdR_i, \tdrho_i) = W(R_i,\rho_i).
\eeqn
As in \cite{Friedman:2005:AMM} one can show that
\beq\label{6.7}
\|u_1-u_2\|_{C_{\xi,t}^{2+\alpha, 1+\alpha/2}(\bar G_{\tau})} \leq C \|(R_1-R_2, \rho_1-\rho_2)\|_Y,
\eeq
from which one can easily deduce that
\beq\label{6.8}
\|\frac{d}{dt}(\tdR_1-\tdR_2)\|_{C^{2+\alpha}[0,\tau]} \leq C \|(R_1-R_2, \rho_1-\rho_2)\|_Y,
\eeq
and
\beqn
\|(\tdM_1-\tdM_2, \tdK_1-\tdK_2)\|_{C_{\xi,t}^{2+\alpha, 1+\alpha/2}(\bar G_{\tau})} \leq C \|(R_1-R_2, \rho_1-\rho_2)\|_Y.
\eeqn
Using arguments as in the proof of Lemma \ref{thm4.9} and \ref{4.12} (iii) and noting that
$\tdrho_1-\tdrho_2=0$ at $t=0$, we derive the estimate
\beqn
\|\tdrho_1 - \tdrho_2 \|_{C_{\xi,t}^{2+\alpha, 1+\alpha/2}(\bar G_{\tau})} \leq C \|(R_1-R_2, \rho_1-\rho_2)\|_Y.
\eeqn
Recalling also (\ref{6.8}) and the fact that $\tdR_1-\tdR_2 = 0$ at $t=0$, we deduce, analogously to (\ref{6.6}), that
\beqn
\|(\tdR_1-\tdR_2, \tdrho_1-\tdrho_2)\|_Y \leq C \tau^{1/2}\|(R_1-R_2, \rho_1-\rho_2)\|_Y.
\eeqn
Hence if $\tau$ is sufficiently small then $W$ is a contraction. We have thus established existence and uniqueness for a small time interval $0\leq t\leq \tau$.

In order to prove existence and uniqueness for all $t>0$ we suppose that such a global solution does not exist and derive a contradiction. Suppose that a unique solution exists for $0\leq t<T$ but not for a larger time interval. We then use the \textit{a priori} estimates of Lemma \ref{4.12} combined with local existence and uniqueness to extend the solution to a larger interval $0\leq t< T+\tau$, which is a contradiction. \q

\section{Ischemic wounds do not heal}\label{7}

In this section we prove that if the parameter $\gm$ in the oxygen equation (\ref{2.4}) and the boundary conditions (\ref{2.13}) -- (\ref{2.18})  is near 1 then $R(t)=\mbox{const.}>0 $ for all $t$ sufficiently large, that is, ischemic wounds do not heal.

For any function $u(r,t)$ we introduce the integral
\beq\label{7.1}
I_u(t) =\int_{R(t)}^L r u(r,t)dr.
\eeq
Using (\ref{2.11}), (\ref{2.12}) we obtain
\beaa
\dd{}{t}\left(\int_{R(t)}^L r u(r,t)dr \right) &=&  \int_{R(t)}^L r \pd{u(r,t)}{t}dr - R(t) u(R(t),t) \dot{R}(t)\\
&=& \int_{R(t)}^L r \pd{u(r,t)}{t}dr + Lu(L,t)v(L) - R(t) u(R(t),t) v(R(t))\nonumber \\
&=& \int_{R(t)}^L r \pd{u}{t} dr + \int_{R(t)}^L  \pd{}{r}(ru  v) dr, \nonumber
\eeaa
or
\beq\label{7.2}
\dd{}{t}I_u(t) \int_{R(t)}^L r \left(\pd{u}{t}   + \frac{1}{r}\pd{}{r}(ru  v) \right) dr.
\eeq
This formula will be used in subsequent lemmas.

For clarity we shall denote the solution $u$ by $u_{\gm}$, and consider first the case $\gm = 1$.

\begin{lma}\label{thm7.1}
There holds:
\beq\label{7.3}
I_{w_1}(t)=\int_{R_1(t)}^L  r w_1(r,t)dr \leq C e^{-\lambda_{wm}t}, \quad C= I_{w_1}(0).
\eeq
\end{lma}
\pf
Multiplying Equation (\ref{2.4}) by $r$ and integrating over  $r\in (R_{\gm}(t), L)$,
we obtain,
\beaa
&& \hspace*{0cm}\dd{}{t}\left(\int_{R_{\gamma}(t)}^L r w_{\gamma}(r,t)dr \right) = LD_w \pd{w_{\gamma}}{r}(L) -R(t)D_w \pd{w_{\gamma}}{r}(R(t)) \\
&& +\int_{R(t)}^L r\left\{
k_{w} b_{\gamma} \big( (1-\gm)w_b - w_{\gamma} \big) - \left[ \big(\ld{wf} f_{\gamma}+ \ld{wm} m_{\gamma}\big)  \left(1 + \dfrac{\ld{ww} p_{\gamma}}{1+ p_{\gamma}}\right) + \ld{wm} \right] w_{\gamma}
\right\}dr
\eeaa
so that, for $\gamma = 1$,
\beaa
\dd{}{t} I_{w_1}(t)\leq -\ld{wm}I_{w_1}(t),
\eeaa
and \re{7.3} follows. \q

\begin{lma}\label{thm7.2}
There holds:
\beqn
I_{f_1}(t)=\int_{R_1(t)}^L  r f_1(r,t)dr  \rightarrow 0 \as t\to\infty.
\eeqn
\end{lma}
\pf
Multiplying Equation (\ref{2.8}) with $\gm=1$ by $r$ and integrating over  $r\in (R_1(t), L)$ we obtain, after using the boundary conditions (\ref{2.16}) and (\ref{2.23}),
\beaa
\dd{}{t}I_{f_1}(t)
&=&  \int_{R_1(t)}^L r\left\{
k_f G_f(w_1)   f_1 \left(1-\frac{f_1}{f_m}\right) - \ld{f} f_1 (1 + \lambda_d D(w_1) )
\right\}dr \\
&\leq& C I_{w_1} (t) - \ld{f} I_{f_1}(t).
\eeaa
Recalling (\ref{7.3}) we deduce
\beq
I_{f_1}(t) \leq (C_1 t+C_2) e^{-\min\{\ld{wm},\ld{f}\}t}\to 0 \as t\to \infty. \q
\eeq

\begin{lma}\label{thm7.3}
There holds:
\beqn
I_{\rho_1}(t)=\int_{R_1(t)}^L  r \rho_1(r,t)dr \rightarrow 0 \as t\to\infty.
\eeqn
\end{lma}
\pf As in the proof of Lemma \ref{thm7.2} one can easily derive the inequality
\beq
\label{7.5}
I_{\rho_1}(t) \leq (C_1 t+C_2) e^{-\min\{\ld{wm},\ld{\rho}\}t} \rightarrow 0 \as t\to\infty. \q
\eeq

From the definition of $Q(r)$ in \re{3.1} and Lemma \ref{thm7.3} we obtain:
\begin{lma}\label{7.4}There holds:
\beqn
Q_1(t) = I_{P_1}(t)=\int_{R_1(t)}^L  r P_1(r,t)dr \rightarrow 0 \as t\to\infty.
\eeqn
\end{lma}

\medskip
We next prove:

\begin{lma}\label{thm7.5}
There exists a constant $C$ such that
\beqn
\max_{R_1(t)\leq r\leq L} w_1(r,t) \leq C e^{-\ld{wm}t/2} \quad \mbox{for all}\;\; t>0.
\eeqn
\end{lma}
\pf
For $\gamma=1$, the oxygen equation can be written in the form
\beqn
 \pd{w_1}{t} - \divp{D_w \pd{w_1}{r} } + v \pd{w_1}{r} +S_1(r,t)w_1 = 0,
\eeqn
where
\beqn
S_1(r,t)=\left[ k_{w} b_1 + \big(\ld{wf} f_1+ \ld{wm} m_1\big)  \left(1 + \dfrac{\ld{ww} p_1}{1+ p_1}\right) + \ld{wm} + P_1(r,t) - \frac{2Q_1(t)}{L^2+R_1(t)^2} \right] .
\eeqn
By Lemma \ref{7.4}, there exists a $t_1$, such that, when $t\geq t_1$, $ 2Q_1(t)/(L^2+R_1(t)^2) \leq  \ld{wm}/2$.
Hence $$S_1(r,t)\geq \ld{wm}/2,$$
and by the comparison lemma,
\beqn
w_1(r,t) \leq \max_{R_1(t)\leq r\leq L} w_1(r,t_1) e^{-\ld{wm}t/2 }.\q
\eeqn

\begin{lma}\label{7.6}There exists a positive constant  $F^{\ast}_1$, $F^{\ast}_1 \geq f_m$, such that
\beqn
  f_1 \leq F^{\ast}_1 \quad \mbox{for all}\;\; R_1(t)\leq r\leq L, \; t>0.
\eeqn
\end{lma}
\pf From Lemmas \ref{7.4} and \ref{thm7.5} it follows that there exists a $t_1>0$ such that,  for all $t\geq t_1$,
$$\frac{2}{L^2+R_1(t)^2}Q_1(t)+k_f G_f(w_1)(1-\frac{f_1}{f_m})\leq\ld{f}/2.$$
Using this in (\ref{2.8}) and setting
\beqn
\bar{f_1} = \max_{0\leq t\leq t_1, R_1(t)\leq r\leq L} f_1(r,t),
\eeqn
we deduce by the comparison lemma that
\beqn
f_1(r,t)\leq \max\{\bar{f}_1, f_m\}\;\;\;\; \mbox{for all}\;\; t\geq t_1. \q
\eeqn

We next improve Lemma \ref{thm7.3}:
\begin{lma}\label{7.7}
\beqn
\max_{R_1(t)\leq r\leq L} \rho_1(r,t) \rightarrow 0 \as t\to\infty.
\eeqn
\end{lma}
\pf By Lemma \ref{7.4}
\beqn
\frac{2}{L^2+R_1(t)^2}Q_1(t)\leq\ld{\rho}/2 \iif t\geq t_1.
\eeqn
Using also Lemmas \ref{thm7.5} and \ref{7.6} we obtain
\beqn
\frac{D\rho_1}{Dt} \leq - \frac{\ld{\rho}}{2}\rho_1 + \frac{F_1^{\ast}k_{\rho}}{K_{w\rho}}
\max_{R_1(t)\leq r\leq L} w_1(r,t) e^{-\ld{wm}t/2} \quad\mbox{for all}\;\; t\geq t_1,
\eeqn
where $D/Dt$ is the derivative along the characteristic curves, and assertion of the lemma follows. \q

Lemma \ref{7.7} implies that $P_1\equiv 0$ for all $t$ sufficiently large, say, for $t\geq T_1^{\ast}$. Hence also $Q_1(t) \equiv 0$ if $t\geq T_1^{\ast}$. Recalling (\ref{3.6}) we conclude:

\begin{lma}\label{7.8}
There exists $R_1^{\ast}>0$ and $T_1^{\ast} > 0$ such that
\beaa
&&R_1(t) > R_1^{\ast} \qquad \mbox{for all}\;\; 0\leq t < T_{1}^{*},  \\
&&R_1(t) \equiv R_1^{\ast} \qquad \mbox{for all}\;\; t\geq T_{1}^{*}.
\eeaa
\end{lma}

We  next extend this result to all $\gm$ near 1.

\begin{thm}\label{7.9}
For any $0\leq 1-\gamma \ll 1$, there exists $R_{\gm}^{\ast}>0$ and $T_{\gamma}^{\ast} > 0$ such that
\beaa
&&R_{\gamma}(t) > R_{\gamma}^{\ast} \qquad \mbox{for all}\;\; 0\leq t < T_{{\gamma}}^{*}, \\
&&R_{\gamma}(t) \equiv R_{\gamma}^{\ast} \qquad \mbox{for all}\;\; t\geq T_{{\gamma}}^{*}.
\eeaa
\end{thm}
\pf Since the estimates of Lemma \ref{4.12} hold uniformly in $\gm$,  any sequence $\gm_i\to 1$ has a subsequence for which the solution $u_{\gm}$ of (\ref{2.2}) -- (\ref{2.28}) converges in $\Omega_{\tau}$, for any $\tau>0$,  to a solution $u_1$ of (\ref{2.2}) -- (\ref{2.28}) with $\gm=1$; the convergence is in the norms of Lemma (\ref{4.12}) with $\alpha$ replaced by any $0<\beta<\alpha$. Since (by Theorem \ref{6.1}) the solution of  (\ref{2.2}) -- (\ref{2.28}) with $\gm = 1$ is unique, we conclude that as $\gm \to 1$ the solution $u_{\gm}$ converges to $u_1$. It follows that
\beqn
\rho_{\gamma}(r,\bar{t}_1 ) \leq \frac{3}{4}, \quad w_{\gm} (r, \bar{t}_1) < \eta_0, \quad f_{\gamma}(r, \bar{t}_1 ) \leq F_1^{\ast} + 1, \quad R_{\gamma}(\bar{t}_1 ) \geq  R_{1}^{\ast}/2
\eeqn
if $\bar{t}_1$ is large enough, provided $\gm\in (\gm_0,1)$ and $1-\gm_0$ is small enough; here $\eta_0$ is chosen small enough so that
\beq\label{thm7.6}
\frac{2 \eta_0 k_{\rho}(F_1^{\ast} + 1)}{K_{w\rho}} \leq \frac{3}{4}\ld{\rho},
\eeq

Let $[\bar{t}_1 , t_{\gamma} )$ be the maximal interval such that
\beqn
\rho_{\gamma}(r,t)<1, \quad \forall t \in [\bar{t}_1 , t_{\gamma} ),
\eeqn
We want to prove that $t_{\gamma}  = +\infty$.
Noting that $Q_{\gamma}(t) \equiv 0$ for  $\bar{t}_{1} \leq t < t_{\gamma}$, we also have $v_{\gamma}(r,t)\equiv 0$ and $R_{\gm}(t)\equiv R_{\gm}(\bar{t}_1)$ for $\bar{t}_1 < t < t_{\gm}$.

Let $W(r,t) = \eta_1 (r-\bar{R})^2 + \eta_0$ where $\bar{R}=R(\bar{t}_1)$, $\eta_1 = (1-\gm)/\bar{A}$ and $\bar{A} = 2\gm L (L- \bar{R})$. Then $(\partial W/\partial r)(\bar{R}, t) = 0$ and
\beqn
(1-\gm)(W-1)+\gm L \pd{W}{r} > 0 \at r=L
\eeqn
if $1-\gm$ is small enough. Also
\beqn
W_t - D_w \Delta W \geq k_w b \big((1-\gm)w_b - W\big)-\ld{wm} W \iif \eta_1 \ll \eta_0,
\eeqn
that is, if $\gm$ is restricted to a very small subinterval $(\gm_1,1)$ of $(\gm_0,1)$. By the comparison lemma we then get

\beqn
w_{\gamma}(r,t)\leq W(r,t)   \quad \mbox{for}\;\; t \in [\bar{t}_1 , t_{\gamma} ),
\eeqn
and, in particular,
\beq\label{thm7.7}
w_{\gamma}(r,t)\leq 2\eta_0    \quad \mbox{for}\;\;  t \in [\bar{t}_1 , t_{\gamma} ).
\eeq

From (\ref{2.2}), \re{thm7.6} and \re{thm7.7} we then obtain, for $\gamma \in (\gamma_1, 1)$,
\beqn
\frac{D\rho}{Dt} \leq \ld{\rho}\left
(\frac{3}{4}-\rho\right) \quad \mbox{for}\;\; t \in [\bar{t}_1 , t_{\gamma} ),
\eeqn
so that
\beqn
\rho_{\gamma}(r,t)\leq \frac{3}{4}, \quad \mbox{for}\;\; t \in [\bar{t}_1 , t_{\gamma} ).
\eeqn
This implies that $t_{\gamma} = +\infty$, and consequently $Q_{\gm}(t)=0$ for all $t>\bar{t}_1$, and the theorem follows. \q

\section{Wounds that do not heal}\label{8}
A wound may be considered to be (completely) healed if $R(t)\to 0$ as $t\to\infty$. Indeed, biologically, if $R(t)$ becomes smaller than, say, $10\; \mu$m (which is roughly the diameter of a cell), no cell can move in to occupy the remaining open space of the wound.
We say that a wound does not heal if
\beq
\lim_{t\to\infty} R_{\gm}(t) = R_{\gm}^{\ast} >0. \label{8.1}
\eeq
In Section \ref{7} we proved that if $\gm$ is near 1 then the wound does not heal and, moreover, $R_{\gm}(t)$ becomes constant for all $t$ large enough.
In this section we want to explore some of the implications of \re{8.1}. In particular we show that in wounds that do not heal, the concentration of oxygen and the density of ECM cannot exceed those of a healthy tissue as $t\to \infty$.

\begin{thm}\label{thm8.1}
If (\ref{8.1}) holds then
\beq
\limsup_{t\to\infty} f_{\gm}(r,t)\leq f_m, \label{8.2}
\eeq
\beq
\limsup_{t\to\infty} w_{\gm}(r,t)\leq \max\{1,(1-\gm) w_b\}, \label{8.3}
\eeq
\beq
\lim_{t\to\infty} \esssup  \rho_{\gm} (r,t) \leq 1. \label{8.4}
\eeq
\end{thm}
\pf
By (\ref{3.6}) and (\ref{3.1}), the function $Q_{\rm}(t)$ satisfies:
\beq\label{8.5}
Q_{\gm}(t) = - \frac{L^2 + R_{\gm}^2(t)}{2R_{\gm}(t)} \dot{R}_{\gm}(t).
\eeq
Integrating over $(0, \infty)$ and recalling (\ref{8.1}), we conclude that
\beqn
\int_0^{\infty} Q_{\gm}(t)dt  = \int_{R_{\gm}^{\ast}}^{R(0)} \frac{L^2 + z^2}{2z} \,dz = \frac{L^2}{2} \log\left(\frac{R(0)}{R_{\gm}^{\ast}}\right) + \frac{R^2(0) - (R_{\gm}^{\ast})^2 }{4}< \infty.
\eeqn

We next prove
\beq
f(r,t) \leq C, \qquad \mbox{for} \; R_{\gm}(t)\leq r\leq L, \;0<t<\infty. \label{8.6}
\eeq
By (\ref{3.8}) we can rewrite the left-hand side of (\ref{2.8}) in the form
\beq\label{8.7}
\pd{f_{\gm}}{t} + v_{\gm} \pd{f_{\gm}}{r} + f_{\gm}\left( P_{\gm}(r,t) - \frac{2}{L^2 + R_{\gm}^2(t)} . Q_{\gm}(t)\right).
\eeq
Hence the function $g(t) = f_m e^{\int_0^t \frac{2}{L^2}Q_{\gm}(s)ds}$ is a supersolution of (\ref{2.8}) and, by the comparison lemma,
$$ f_{\gm}(r,t) \leq g(t), \quad R(t)\leq r\leq L, \;\; t>0.$$
Since, by (\ref{8.5}), $g(t)$ is uniformly bounded, (\ref{8.6}) follows.

We next prove that
\beq
|\dot{Q}_{\gm}(t)| \leq C  \quad \mbox{for all  } t>0. \label{8.8}
\eeq
We write (\ref{2.2}) in the form
\beqn
\pd{ (\rho_{\gm}-1) }{t} + \frac{1}{r} \pd{}{r} \Big(r(\rho_{\gm} - 1) v_{\gm} \Big) = -\frac{1}{r}\pd{(rv_{\gm})}{r} + \frac{k_{\rho} w_{\gm} }{w_{\gm}+K_{w\rho}} f_{\gm} (1-\frac{\rho_{\gm}}{\rho_m}) - \lambda_{\rho} \rho_{\gm}  \triangleq M_{\gm},
\eeqn
or
\beqn
r\pd{P_{\gm}}{t} + \pd{}{r} (rP_{\gm}v_{\gm}) = r\beta M_{\gm} I_{\{ (r,t): \;\rho_{\gm}(r,t)>1 \}}.
\eeqn

By \re{3.8} (or \re{thm4.7}), \re{8.6} and the bound $\rho_{\gm} \leq \rho_m$, we see that the right-hand side is uniformly bounded in $(r,t)$.
Hence, by integration over $R(t)\leq r\leq L$,
\beq\label{8.9}
\int_{R(t)}^L r \pd{P_{\gm(r,t)}}{t} dr \mbox{   is uniformly bounded.}
\eeq
Next, by the definition of $Q_{\gm}(t)$ in \re{3.1},
\beqn
\dot{Q}_{\gm}(t) =  - \dot{R}(t) R(t) P(R(t),t) + \int_{R(t)}^L r \pd{P_{\gm(r,t)}}{t} dr,
\eeqn
and hence, upon using \re{8.9} and the uniform boundedness of $\dot{R}(t)$, the assertion \re{8.8} follows.

From (\ref{3.6}) and (\ref{8.8}), we obtain the estimate
\beq
|\ddot{R}_{\gm}(t)|\leq C. \label{8.10}
\eeq

Using the interpolation estimate (see \cite{LiebermanPDE}, Page 48)
\beqn
\| \dot{R}_{\gm} \|_{C^{\alpha}[t^{\ast}, t^{\ast}+1]} \leq C \| R_{\gm}- R_{\gm}^{\ast}  \|^{\frac{1+\alpha}{2}} _{W^{2,\infty}[t^{\ast}, t^{\ast}+1]} \cdot
\| R_{\gm} - R_{\gm}^{\ast} \|^{\frac{1-\alpha}{2}} _{L^{\infty}[t^{\ast}, t^{\ast}+1]} \rightarrow 0
\eeqn
and noting that the last factor converges to zero as $t^{\ast}\to 0$, we obtain
\beq
\lim_{t^{\ast}\to\infty} \| \dot{R}_{\gm} \|_{C^{\alpha}[t^{\ast}, t^{\ast}+1]} = 0,
\eeq
and then, by (\ref{3.6}), also
\beq\label{8.12}
\lim_{t^{\ast}\to\infty} \| Q_{\gm} \|_{C^{\alpha}[t^{\ast}, t^{\ast}+1]} = 0, \qquad \forall\, 0<\alpha<1.
\eeq

From (\ref{8.12}) and (\ref{8.5}) it easily follows that
\beq\label{8.13}
Q_{\gm}(t)\to 0 \as t\to \infty,
\eeq
hence there exists a $T>0$ such that
\beq
\frac{2}{L^2} Q_{\gm}(t) \leq \frac{\ld{f} }{2} \quad \mbox{if}\;\; t>T.
\eeq
Writing the left-hand side of \re{2.8} in the form \re{8.7} and using \re{8.13}, we
can then apply the comparison lemma to $f_{\gm}$   to conclude that
\beqn
f_{\gm}(r,t) \leq f_m + \max_{R_{\gm}(T)\leq r\leq L} f_{\gm}(r,T) e^{-\frac{\ld f}{2}(t-T)},
\eeqn
and hence (\ref{8.2}) follows.

Similarly one can prove, by comparison, the estimate (\ref{8.3}). Finally, (\ref{8.4}) follows from \re{8.13}. \q

\section{Simulations and a conjecture}\label{9}

\begin{figure}[ht]
\includegraphics[width=4in]{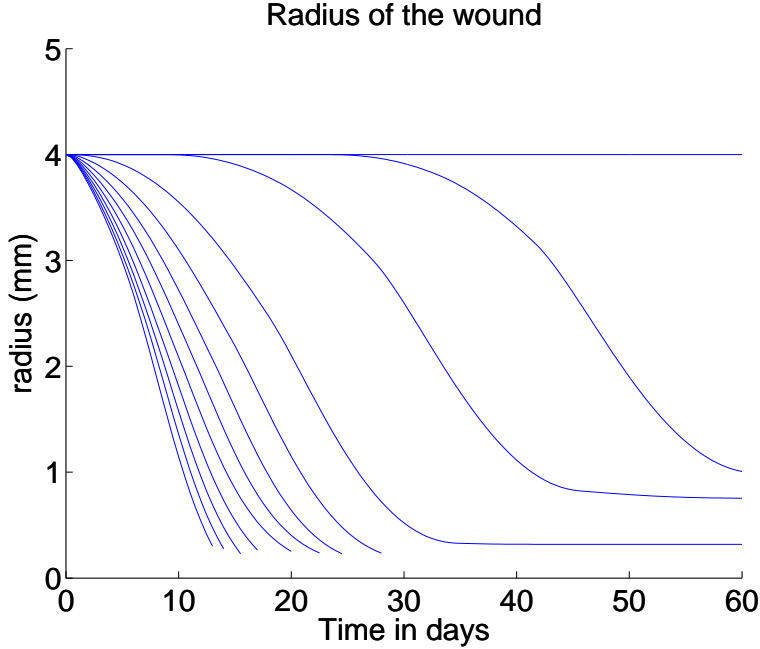}
\caption{The radius of the wound as a function of time for different values of $\gm$. From left to right: $\gm = 0, 0.1, 0.2, \ldots, 0.8, 0.9, 0.92, 0.95, 1$. Other parameters used are the same as in \cite{Xue:2009:MMI}; the nondimensionalized values are: $L=5$, $R_0=8/3$, $\rho_m=2$, $K_{w\rho}=K_{wf}0.25$, $k_{\rho}=5/16$, $\ld{\rho}=0.1$, $\beta=10$, $D_w=0.5$, $D_p=D_e=1$, $D_m=D_f=5\times 10^{-2}$, $D_n=10^{-3}$, $D_b=7\times 10^{-4}$, $\chi_m=\chi_f=0.1$, $\chi_n=1$, $m_m=f_m=n_m=10$, $A=0.1$, $w_b=2$, $k_w=4.39$, $\ld{wf}=0.227$, $\ld{wm}=4.16$, $\ld{d}=2$, $k_f=5.78\times 10^{-3}$, $\ld{f}=5.2\times 10^{-3}$, $k_{nb}=k_n=k_b/10=2.16\times 10^{-2}$, $\ld{nn}=100\ld{nb}=2.25$, $k_{sg} = 6.25\times 10^{-2}$.}\label{9.1}
\end{figure}

We simulated the radius $R_{\gm}(t)$ of the wound for different values of $\gm$ using the nondimensional
parameters of the system \re{2.2} -- \re{2.28} that were chosen on the basis of experimental results \cite{Xue:2009:MMI}.
In Figure \ref{9.1} we present simulation results in the original dimensional variables with $L = 7.5$ mm and initial wound radius $R_0 = 4$ mm.
The computation was manually stopped when the wound
became 98\% closed. From the figure we see that as $\gm$ increases, the wound closes slower, and when $\gm$ is close
to 1, the wound radius stops decreasing after a certain time.

We conjecture that if the parameters of the system (\ref{2.2}) -- (\ref{2.28}) are chosen on the basis of experimental results, as in \cite{Xue:2009:MMI}, then there exists a parameter value $\gm_{\ast}$ such that (\ref{8.1}) holds if $\gm_{\ast}<\gm\leq 1$ and
\beqn
\lim_{t\to\infty} R_{\gm}(t) = 0 \iif \gm<\gm_{\ast}.
\eeqn
If this conjecture is true then, in particular,
\beqn
\lim_{t\to\infty} R_{\gm}(t) = 0 \iif \gm=0.
\eeqn
But even this assertion is still an open question. We can only prove, for the system \re{2.2} -- \re{2.28}, with general parameters, the following result.

\begin{thm}
 If $\gamma =0 $, then
\bea
  && \rho(L,t)>1, \mm  0<t<\infty, \label{thm9.1}\\
  && \dot{R}(t) < 0, \mm 0<t<\infty, \label{9.2}\\
  && Q(t) >0, \mm 0<t<\infty.\label{9.3}
\eea
\end{thm}
\pf Using the boundary conditions $w(L,t)=1, \; f(L,t)=1$, $v(L,t)=0$
and \re{2.26}, we obtain from \re{2.2} at $r=L$ the relation
\beaa
\pd{\rho(L,t)}{t} + \rho(L,t)\pd{v}{r}(L,t) &=& \frac{k_{\rho} }{1+K_{w\rho}} \left(1-\frac{\rho}{\rho_m}\right) - \frac{k_{\rho} }{1+K_{w\rho}} \left(1-\frac{1}{\rho_m}\right) \rho  \\
&=& -\frac{k_{\rho} }{1+K_{w\rho}}(\rho-1),
\eeaa
and from \re{3.8},
\beqn
\pd{v}{r}(L,t) = P(L,t) - \frac{2}{L^2 + R^2(t)} Q(t).
\eeqn
Hence,
\be\label{9.4}
  \frac{\partial\rho(L,t)}{\partial t} = - c_0 \Big(\rho(L,t)-1\Big)
 -\beta\rho(L,t)\Big(\rho(L,t)-1\Big)_+ +
\frac{2\rho(L,t)}{L^2+R(t)^2} Q(t).
\ee
where $c_0$ is a positive constant.

Using the initial conditions
\beaa
&& \rho(r,0)\equiv 1, \m
w(r,0)\equiv 1, \m f(r,0)\equiv 1,
\m m(r,0)\equiv 0, \m v(r,0)\equiv 0 \m\text{for } R_0<r<L,\\
&&
b(r,0)\equiv 1, \m p(r, 0)\equiv 0 \m \text{for }   R_0+\eps_0\le r <L,
\eeaa
in \re{2.2} and \re{2.4} and recalling the relations (\ref{2.26}) and (\ref{2.27}), we find that
\bea
 && \frac{\partial\rho(r,0)}{\partial t}\equiv 0, \mm R_0<r<L, \label{9.5}\\
 && \frac{\partial w(r,0)}{\partial t}\equiv 0, \mm R_0+\eps_0\le r<L.\label{9.6}
\eea

Using \re{2.25} we also obtain (upon recalling \re{2.28}) that
\be
\frac{\partial f(r,0)}{\partial t} = \frac{\chi_f}{r}H(1-\frac{1}{f_m}) p_0'(r) \Big/
 \sqrt{1+k_{sg}\Big(\frac{k_{pb}}{D_p}\Big)^2}> 0, \mm R_0<r<R_0+\eps_0.
\label{9.7}
\ee

Differentiating Equation (\ref{2.2}) in $t$ and using (\ref{9.5}) -- (\ref{9.7}) and the $C_{r,t}^{2+\alpha, 1+\alpha/2}$ regularity of $w$,  we deduce that
\beqn
\frac{\partial^2 \rho(r_0,0)}{\partial t^2} >0 \quad \mbox{for} \;\; 0< R_0+\eps_0 -r_0 \ll 1
\eeqn
This implies that , for $ 0< R_0+\eps_0 -r_0 \ll 1$ and $0<t\ll 1$,
\[
  \rho(r_0,t)>1,
\]
and hence
\be\label{9.8}
 Q(t)>0 \m\text{for } 0<t\ll 1.
\ee
Since $\rho(0,L)=1$, from \re{9.4} and \re{9.5} it follows that
\be\label{9.9}
 \rho(0,t)>1
\ee
for {\em all} $0<t<\infty$. This in turn implies that $Q(t)>1$ for {\em all} $0<t<\infty$, hence $\dot{R}(t)<0$ and (by (\ref{9.4})) $\rho(L,t)>1$ for {\em all} $0<t<\infty$.\q

\section{Conclusion}\label{10}
In this paper we established existence and uniqueness of a solution to a free boundary problem which
models ischemic wound healing. The ischemic condition is described in terms of a parameter $\gm$
($0 \leq \gm \leq 1$) which appears as a coefficient in a Robin boundary condition for the various cells and chemical
densities. We also proved that under extreme ischemic conditions ($\gm$ near 1) the open wound stops decreasing in finite time.
When the parameters of the system are taken on the basis of biological experiments, simulations show that there is a parameter
$\gm_{\ast}$ such that the wound heals if $0\leq \gm < \gm_{\ast}$ and does not heal if $\gm_{\ast}<\gm\leq 1$. This assertion remains a
challenging mathematical open problem. Future work should include the introduction of pressure and diabetic conditions in ischemic wounds, as well as inflammatory conditions.

\medskip

{\bf Acknowledgment.} This work was partially supported by National
Science Foundation upon agreement No. 0635561 and the National
Institute of Health (OD) Award
UL1RR025755.


\end{document}